\documentclass[12pt]{article}
\usepackage{amsmath,amsfonts,amsthm,amssymb,indentfirst}

\usepackage{amsmath}
\usepackage{amsthm}
\usepackage{amsfonts}
\usepackage{amssymb}

\input{xy}
\xyoption{all}
\usepackage{amscd}

\theoremstyle{plain}
\newtheorem{lemma}{Lemma}[section]
\newtheorem{cor}[lemma]{Corollary}
\newtheorem{thm}[lemma]{Theorem}
\newtheorem*{thmnonum}{Theorem}
\newtheorem{prop}[lemma]{Proposition}

\theoremstyle{definition}
\newtheorem{defn}[lemma]{Definition}
\newtheorem{example}[lemma]{Example}
\newtheorem{q}[lemma]{Question}

\newtheorem{remark}[lemma]{Remark}

\newcommand{\ann}{\operatorname{\ensuremath {\rm ann}}}
\newcommand{\Max}{\operatorname{\ensuremath {\rm Max}}}

\newcommand{\xym}{\ensuremath \xymatrix@1}

\newcommand{\End}{\operatorname{\ensuremath {\rm End}}}

\newcommand{\Hom}{\ensuremath {\rm Hom}}

\newcommand{\tr}{\operatorname{\ensuremath {\rm tr}}}

\newcommand{\Q}{\ensuremath \mathbb{Q}}
\newcommand{\QQ}{\ensuremath \mathbb{Q}}

\newcommand{\CC}{\ensuremath \mathbb{C}}

\newcommand{\NST}{{\ensuremath {\rm Nil_{\ast}}}}

\newcommand{\C}{\ensuremath {C}}
\newcommand{\Z}{\operatorname{\ensuremath {\it {Z}}}}
\newcommand{\F}{\ensuremath \mathbb{F}}

\newcommand{\HH}{\ensuremath {\mathbb{H}}}
\newcommand{\M}{\ensuremath \mathbb{M}}

\newcommand{\T}{\ensuremath \mathbb{T}}

\newcommand{\ta}{{$($P$)$}}
\newcommand{\tb}{{$($PI$)$}}
\newcommand{\tc}{{$($SR$)$}}
\newcommand{\td}{{$($SI$)$}}
\newcommand{\te}{{$($N$)$}}
\newcommand{\tdom}{{$($D$)$}}
\newcommand{\tred}{{$($R$)$}}
\newcommand{\tnot}{{$($N$)$}}
\newcommand{\tta}{{$($P$)$\ }}
\newcommand{\ttb}{{$($PI$)$\ }}
\newcommand{\ttc}{{$($SR$)$\ }}
\newcommand{\ttd}{{$($SI$)$\ }}
\newcommand{\tte}{{$($N$)$\ }}
\newcommand{\ttdom}{{$($D$)$\ }}
\newcommand{\ttred}{{$($R$)$\ }}
\newcommand{\ttnot}{{$($N$)$\ }}

\begin{document}
\bibliographystyle{alpha}

\title{On minimal extensions of rings}
\author{Thomas J. Dorsey and Zachary Mesyan}
\maketitle

\begin{abstract}  
Given two rings $R \subseteq S$, $S$ is said to be a minimal ring
extension of $R$ if $R$ is a maximal subring of $S$.  In this article,
we study minimal extensions of an arbitrary ring $R$, with particular focus
on those possessing nonzero ideals that intersect $R$ trivially.
We will also classify the minimal ring extensions of prime rings,
generalizing results of Dobbs, Dobbs \& Shapiro, and Ferrand \&
Olivier on commutative minimal extensions.
\end{abstract}

\section{Introduction}

Throughout, all rings are associative with unity $1$, which is preserved by homomorphisms and inherited by subrings.  
Rings which do not necessarily have a unity element will be referred
to as {\emph{rngs}}.  A ring $S$ is said to be a {\emph{ring
    extension}} of a ring $R$ if $R$ is a subring of $S$; in particular, $R$ and $S$ must share the same unity element.  Moreover, 
we will say that $S$ is a {\emph{minimal ring extension}} (or {\emph{minimal extension}}, for short) of $R$ if $R$ is a maximal 
subring of $S$.  Explicitly, this holds whenever there are no subrings
strictly between $R$ and $S$.  

Minimal ring extensions 
have been studied in a number of
papers (a great number of which restrict entirely to the category of
commutative rings) and we will provide a brief summary of some of 
that work.  
Ferrand and Olivier classified the minimal
{\emph{commutative}} extensions of fields in their 1970 paper \cite{F&O}.  Much later, in \cite{D&S1}, Dobbs and
Shapiro classified the minimal {\emph{commutative}} extensions of
integral domains.  In \cite{SSY}, Sato, Sugatani, and Yoshida showed
that for any domain $R$ which is not equal to its quotient field $Q(R)$, each
domain minimal extension of $R$ is an {\emph{overring}} in the sense
that it embeds in $Q(R)$.  In \cite{D&S2}, Dobbs and Shapiro examined
the {\emph{commutative}} minimal extensions of certain non-domains, 
and (in a certain sense) reduced the study of commutative minimal extensions
to extensions of reduced rings.  Other aspects of
commutative minimal extensions are studied in 
\cite{Pic1}, \cite{PicPic}, \cite{Old}, \cite{ayache}, \cite{Dobbs},
\cite{Dobbs2}, \cite{Papick}, \cite{Ouk1}, \cite{Ouk2}, \cite{Ouk3},
\cite{Ouk4}, as well as others.  
Minimal extensions of arbitrary rings, as well as minimal
noncommutative extensions of commutative rings, have received
considerably less attention than their commutative counterparts.
Moreover, most of the papers that do study noncommutative minimal
extensions actually study rings which have a maximal
subring of a prescribed type (e.g., having a certain finiteness
property).  The goal there is generally to show that this implies a finiteness
condition on the larger ring.  The main result in this area, found by
A. Klein in \cite{Klein} and T. Laffey in \cite{Laffey},
independently, is that a ring with a finite maximal subring must be,
itself, finite.  An analogue of this is found in \cite{MaxAlg},
where the authors show that if $k$ is a field, then a $k$-algebra
which has a finite-dimensional maximal subalgebra must be, itself,
finite-dimensional.  Other papers in this area include
\cite{bg} and \cite{BK}.  
On a related topic, extensions of rings minimal among a given class of
rings are studied in \cite{BPR1}, \cite{BPR2}, and \cite{BPR3}, as
well as in other papers by the same authors.  

The main subject of this article is the study of minimal
{\emph{ideal}} extensions $S$ of a ring $R$, namely those minimal
extensions possessing a nonzero ideal $I$ which intersects $R$
trivially.  More specifically, given a ring $R$ and an $R$-rng $I$
(i.e., a rng possessing a compatible $(R,R)$-bimodule structure; 
see Definition~\ref{defineeir}), the ideal extension of $R$ by $I$, 
denoted by $E(R,I)$, is the ring whose underlying abelian group is $R \oplus I$, 
and where multiplication is defined by 
$(r,i) \cdot (r',i') = (rr', ir' + ri' + ii')$, for $r, r' \in R$ and
$i, i' \in I$. 
Ideal extensions are quite common, and a familiar example is that of
an 
``idealization'', which is also called a ``split-null'' or ``trivial'' extension 
(see Section~\ref{EIR}).  The importance of ideal extensions to the
study of minimal extensions is suggested by Proposition~\ref{teaser}
below, which asserts that all non-prime minimal extensions of a prime
ring are ideal extensions (some prime minimal extensions of
prime rings are ideal extensions, as well).  We will describe the
minimal ideal extensions of an arbitrary ring, 
and will use this information to classify all minimal extensions of a
prime ring.  
This work generalizes the classification of {\emph{commutative}}
minimal extensions of domains which was done by Dobbs and Shapiro in
\cite{D&S1} 
(following the earlier classification of commutative minimal extensions of fields performed by Ferrand and Olivier in \cite{F&O}).      

The outline for this article is as follows.  In Sections~\ref{EIR} and \ref{idealtheory}, we
will study the general theory of minimal ideal extensions, among other
things, describing the ideal theory of ideal extensions
(Proposition~\ref{idealdescription}), and using this to characterize
when an ideal extension is (semi)prime
(Propositions~\ref{semiprimeoversemiprime} and \ref{primeidealext}).
We will also find and describe three classes of ideal extensions
which stratify the ideal-theoretic behavior of an ideal extension, and
control whether an ideal extension is (semi)prime.
In Section~\ref{centralsection}, we will examine 
{\emph{central extensions}} (namely, extensions of $R$ which are
generated as left $R$-modules by elements which centralize $R$),
whose behavior closely models that of commutative minimal extensions
of commutative rings, and we will characterize when a minimal ideal
extension is a central extension.  

In Section~\ref{primesection}, we will prove the following classification of the minimal extensions of arbitrary prime rings.  

\begin{thmnonum}  Let $R$ be a prime ring.  Then, up to $R$-isomorphism, every minimal 
extension of $R$ must be one of exactly one of the following five forms.
\begin{enumerate}
	\item[\ta]  A prime minimal extension of $R$, all of whose nonzero ideals intersect $R$ nontrivially. 
	\item[\tb]  $E(R,I)$ for some minimal $R$-rng $I$ such that $\,\Hom_R(I,R) = 0$, $I^2 \ne 0$, and $\, \ann_R(I) = 0$.  
	\item[\tc]  $E(R,I)$ for some minimal $R$-rng $I$ such that
	$\, \Hom_R(I,R/\ann_R(I)) = 0$, $I^2 \ne 0$, and $\, \ann_R(I) \ne 0$.
	\item[\td] $E(R,I)$, where $I$ is a minimal ideal of $R/P$
	for some prime ideal $P$ of $R$.  
	\item[\te] The trivial extension $R \propto M$ for some simple $(R,R)$-bimodule $M$.  
\end{enumerate}
Extensions of the forms \tta and \ttb are prime; those of forms \ttc and
\ttd are semiprime, but not prime; and those of form \tte are not
semiprime.  In each case where
they occur, $I$, $M$, and $P$ are unique, up to $R$-isomorphism,
$(R,R)$-bimodule isomorphism, and equality, respectively.

$($The labels are intended to mean:  \tta = prime; \ttb = prime, ideal
extension; 
\ttc = semiprime, reducible; \ttd = semiprime, subdirectly
irreducible; \tte = not semiprime.$)$     
\end{thmnonum}
Our result generalizes the aforementioned
classification of commutative minimal extensions of integral
domains appearing in \cite[Theorem 2.7]{D&S1}.  
Specifically, in the case of central extensions, the above result 
reduces to the following,
which, for $R$ commutative, is essentially identical to the main result of \cite{D&S1}. 
\begin{thmnonum} Let $R$ be a prime ring.  Then, up to
  $R$-isomorphism, every central minimal extension of $R$ must be of exactly one of the following three forms.
\begin{enumerate}  
  \item [\ta] A prime minimal extension of $R$, all of whose nonzero ideals
  intersect $R$ nontrivially.  
  \item [\td] $R \times R/M$ for some maximal ideal $M$ of $R$.
  \item [\te] $R \propto R/M$ for some maximal ideal $M$ of $R$. 
\end{enumerate} 
The maximal ideal $M$, where it appears, is determined by the
$R$-isomorphism type of the extension.
\end{thmnonum}  
Despite being almost identical in statement to \cite[Theorem
  2.7]{D&S1}, this result was proved by a very different method, since the techniques used in \cite{D&S1} (primarily localization) do not carry over to the noncommutative setting.  As illustrated above, the general classification (Theorem~\ref{primeext})
of minimal extensions of an arbitrary prime ring is quite a bit
different from the central case; in fact, even commutative domains (as
well-behaved as $k[x]$, for a field $k$) can have interesting
noncommutative minimal extensions of a flavor entirely different from
the rings appearing in the Dobbs and Shapiro classification, and the
types which do appear in the central case are degenerations of the
corresponding cases appearing in the general classification.  

Finally, in Section~\ref{simplesection}, we will
classify the minimal extensions of simple rings, and we
will present two examples due to George Bergman of the types of
minimal extensions which cannot exist in the central case.  In
particular, we will produce non-simple prime minimal extensions of
certain fields and noncommutative semiprime non-prime minimal
extensions of commutative domains.  

\subsection*{Acknowledgements}  
The authors would like to thank George Bergman, Alex Diesl, Danny
Goldstein, and Murray Schacher for helpful discussions regarding this
material.

\section{Minimal ideal extensions} \label{EIR}

Given a ring $R$, an $R$-{\em ring} (resp. $R$-{\em rng}) $I$ is a ring
(resp. rng) that is a unital $(R,R)$-bimodule,
for which the actions of $R$ are compatible with multiplication in
$I$.  That is, $r(xy) = (rx)y$, $x(ry) = (xr)y$, and $(xy)r = x(yr)$ for every $r \in R$
and $x,y \in I$.  Note that a ring homomorphism $R \rightarrow I$
equips $I$ with the structure of an $R$-ring in a natural way; in
particular, in this way every ring extension of $R$ may be viewed as an $R$-ring.    

We will call a nonzero $R$-rng {\emph{minimal}} if it has no proper nonzero
$R$-subrngs.  We note that if $I$ is a minimal $R$-rng with $I^2 \ne 0$, then $I$ is simple as a rng
(i.e., it has precisely two ideals).  This can be proved using essentially an argument found in \cite[Lemma 2(i)]{BK}.  
(The annihilators $\{x \in I: Ix = 0\}$ and $\{x \in I: xI = 0\}$ must
each be zero, since $I^2 \ne 0$ implies that each is a proper $R$-subrng of $I$.  Thus, if $J$ is a nonzero ideal of $I$, then $JI \ne 0$,
and hence $IJI \ne 0$.  But $IJI$ is then a nonzero $R$-subrng of $I$, so $IJI = I$,
by minimality.  On the other hand, $J$ is an ideal of $I$, implying that $IJI
\subseteq J$. So we conclude that $J = I$, and hence that $I$ is
simple.)    
  
Given two $R$-rngs $I$ and $J$, $\Hom_R(I,J)$ will denote the set of 
$R$-homomorphisms $\varphi: I \rightarrow J$ (where an $R$-homomorphism is a homomorphism of $R$-rngs that is also an 
$(R,R)$-bimodule homomorphism).  
Given an $R$-rng $I$, 
$\ann(I_R) = \{x \in R: Ix = 0\}$ will denote the right annihilator of $I$ in $R$; 
$\ann({}_R I) = \{x \in R: xI = 0\}$ will denote the left annihilator of $I$ in $R$; 
and we set $\ann_R(I) = \ann(I_R) \cap \ann({}_R I) = \{x \in R: xI = Ix = 0\}$.  Each of these annihilators is a 2-sided ideal of $R$.  

We begin with a basic lemma regarding annihilators.  
\begin{lemma}  \label{minimalann} 
Let $R$ be a ring, and let $I$ be a minimal $R$-rng.  Then $\ann(I_R)$ and $\ann({}_R I)$ are prime $($$2$-sided$)$ ideals of $R$, and hence 
$\ann_R(I)$ is a semiprime ideal.  If $I^2 \ne 0$, then $\ann_R(I) = \ann(I_R) = \ann({}_R I)$, and hence $\ann_R(I)$ is 
prime.  
\end{lemma}
\begin{proof}
Suppose that $A$ and $B$ are ideals of $R$ for which $AB \subseteq \ann(I_R)$, or equivalently, $I(AB) = 0$.  
By minimality, either $IA = 0$, in which case $A \subseteq \ann(I_R)$; or else $IA = I$, in which case 
$0 = I(AB) = (IA)B = I B$, so $B \subseteq \ann(I_R)$.  We conclude that $\ann(I_R)$ is prime.  Similarly, $\ann({}_R I)$
is prime, and hence their intersection $\ann_R(I)$ is semiprime.  

Now, suppose that $I^2 \ne 0$.  By minimality, the $R$-subrng $\ann(I_R) I$ of $I$ is either $0$ or $I$, 
but $(\ann(I_R) I)^2 = 0$, which forces $\ann(I_R) I = 0$, since $I^2 \ne 0$.  We conclude that $\ann(I_R) \subseteq \ann({}_R I)$.  By a similar argument, we conclude that $\ann(I_R) \supseteq \ann({}_R I)$, and hence $\ann(I_R) = \ann( {}_R I) = \ann_R(I)$.     
\end{proof}

Given an $R$-rng $I$, there is a natural way of enlarging $I$ to an
$R$-ring (which is the Dorroh extension when $R = \mathbb{Z}$).  
The construction that follows can be viewed as a functor from the category of $R$-rngs (and
$R$-homomorphisms) to the category of $R$-rings (and
$R$-homomorphisms).  
\begin{defn} \label{defineeir}
Given an $R$-rng $I$, the {\it ideal extension} $E(R,I)$ has the abelian group structure of $R \oplus I$, with
multiplication given by $(r,i) \cdot (r',i') = (rr', ir' + ri' + ii')$.  We identify 
the subring $R \oplus 0$ with $R$, and we identify the $R$-rng $I$ with the ideal $0 \oplus I$ 
of $E(R,I)$.  
It is straightforward to verify that with these operations (and the embedded copy of $R$) 
$E(R,I)$ is an $R$-ring, and it is easy to see that the assignment of $I$ to $E(R,I)$ is functorial.  
\end{defn}

One common instance of this construction is the ``trivial
extension'' (which is also called a ``split-null'' extension or an
``idealization'') 
$R \propto M$, 
of a ring $R$ by an $(R,R)$-bimodule $M$, which is the ring with underlying abelian group structure of $R \oplus M$, 
and multiplication defined by $(r,m) \cdot (r',m') = (rr', rm' + mr')$, where $r,r' \in R$ and $m,m' \in M$.  
Viewing $M$ as an $R$-rng with square zero multiplication, clearly $E(R,M) = R \propto M$.              

As we shall see in Lemma~\ref{produce} below, 
ideal extensions are relevant to the study of minimal extensions in
general.  The next result illustrates that they are truly essential when studying minimal extensions of prime rings.

\begin{prop}  \label{teaser}
Let $R$ be a ring, and let $S$ be a minimal ring extension of $R$ which
has a nontrivial ideal that intersects $R$ trivially.  Then, $S$ is
$R$-isomorphic to an ideal extension of $R$.  In particular, if $R$ is a prime ring and
$S$ is a minimal extension of $R$ which is not prime, then $S$ is
$R$-isomorphic to an ideal extension of $R$. 
\end{prop} 
\begin{proof}
For the first statement, let $I$ be a nonzero ideal of $S$ which
intersects $R$ trivially.  The additive group $R+I$ is a subring of
$S$ which properly contains $R$, so by minimality, $R+I = S$, where
the sum is direct.  It follows easily that $S$ is
$R$-isomorphic to $E(R,I)$.  

Now, suppose that $R$ is prime and $S$ is a minimal extension of $R$
which is not prime.  
Thus, there exist nonzero ideals $I$ and $J$ of $S$ for which $IJ = 0$.  
But then $(R \cap I)(R \cap J) = 0$, so one of the two ideals $R \cap I$ and $R \cap J$ of $R$ 
must be zero.  Without loss of generality, $R \cap I = 0$, and as
above, $S$ is $R$-isomorphic to $E(R,I)$.  
\end{proof}

The next lemma relates the structure of the ideal extension $E(R,I)$ to the structure of the $R$-rng $I$.  
In the case of trivial extensions over commutative rings, this is
simply \cite[Theorem 2.4 and Remark 2.9]{Dobbs}, 
and the proof we give here is similar to the proof appearing there.  
\begin{lemma} \label{produce}  Let $R$ be a ring, and let
  $I$ be an $R$-rng.  The map $K \to E(R,K)$ is a
  one-to-one, inclusion preserving, correspondence between the
  $R$-subrngs of $I$ and the subrings of $E(R,I)$ which
  contain $R$.  Consequently, $E(R,I)$ is a minimal extension of $R$ if and only if $I$ is a minimal $R$-rng.   
\end{lemma}  
\begin{proof}
It is clear that the map sending an $R$-subrng $K$ of $I$ to $E(R,K)
\subseteq E(R,I)$ is inclusion preserving, and sends $R$-subrngs of
$I$ to subrings of $E(R,I)$ which contain $R$.  
The inverse map sends a subring $R \subseteq S \subseteq E(R,I)$ to its image $S^{\sharp}$ under the $(R,R)$-bimodule 
homomorphism projecting $E(R,I)$ to its second component $I$ (this map is not, in general, an $R$-rng homomorphism).  
It is straightforward to see that $S = R \oplus S^{\sharp}$, and from this it follows that the $(R,R)$-subbimodule 
$S^{\sharp}$ of $I$ must be closed under multiplication.  Hence $S^{\sharp}$ is an $R$-subrng of $I$, and therefore the map $K \to E(R,K)$ is a one-to-one correspondence between $R$-subrngs of $I$ and subrings of $E(R,I)$ which contain $R$.    

The final claim is clear.
\end{proof}

\begin{remark}
It follows, as in \cite{Dobbs}, that any ring $R$ has a minimal ring
extension, since for any maximal ideal $M$ of $R$, the
$(R,R)$-bimodule $R/M$, viewed as an $R$-rng $I$ with trivial multiplication, is a
minimal $R$-rng, and hence $E(R,I) \cong R \propto M$ is a minimal ring
extension of $R$.  
\end{remark}

While Proposition~\ref{teaser} gives a strong reason to consider ideal
extensions in the context of the study of minimal extensions, in some
sense minimal ideal extensions must be ``seen'' by semiprime (or
prime) rings.  Indeed, suppose that $E(R,I)$ is a minimal ideal
extension of a ring $R$.  Whether or not $R$ is semiprime, by
Lemma~\ref{minimalann}, $R/\ann_R(I)$
is semiprime (and prime if $I^2 \ne 0$), and $E(R,I)/(\ann_R(I) \oplus
0) \cong E(R/\ann_R(I),I)$.  Moreover, since $\ann_R(I)$ is semiprime
by Lemma~\ref{minimalann}, $\ann_R(I) \supseteq \NST(R)$,
the lower nil (or prime) radical of $R$ (which is the smallest
semiprime ideal of $R$), and $E(R,I)/(\NST(R) \oplus 0)
\cong E(R/\NST(R), I)$, so every minimal ideal extension of $R$ yields
a minimal ideal extension of the maximal semiprime quotient
$R/\NST(R)$.  

In the commutative case, \cite[Theorem 2.1]{D&S2} gives a reduction of
the study of minimal commutative extensions to minimal commutative
extensions of reduced rings.  A fact used in the proof of this theorem is that if $R$ is a subring of a ring $S$, then $R \cap \NST(S) \subseteq \NST(R)$ (\cite[Ex. 10.18A(1)]{fc}). In the noncommutative setting, the situation is
complicated by the fact that  this containment can be strict, whereas, if $R$ lies in the center of $S$, then $R \cap \NST(S) = \NST(R)$.  In particular, a semiprime ring can
have a non-semiprime subring, which cannot occur in the category of
commutative rings.  For an example with minimal ring extensions, let
$k$ be a field, and consider $\T_2(k) \subseteq \M_2(k)$, the subring
of $2 \times 2$ upper triangular matrices in the full ring of $2\times
2$ matrices over $k$.  Comparing
$k$-dimensions shows that this is a minimal ring extension, but
$\M_2(k)$ is semiprime, whereas $\T_2(k)$ is not.  

If $S$ is a minimal extension of a ring $R$ with $\NST(R) = \NST(S)$, then
clearly $S/\NST(S)$ is a minimal extension of $R/\NST(R)$.  
The following lemma, provides an analogue
to the last statement in \cite[Theorem 2.1]{D&S2}, but the two
conditions $\NST(R) = \NST(S)$ and $\NST(S) \not\subseteq R$ do not
exhaust all possibilities, as they do in the commutative case.      
\begin{lemma}  \label{partialreduction} Let $R$ be a ring, and let $S$ be a minimal extension
  of $R$.  If $\, \NST(S) \not\subseteq R$, then $\, \NST(R) = \NST(S) \cap
  R$ and $R/\NST(R) \cong S/\NST(S)$.  
\end{lemma}
\begin{proof}
Suppose that $\NST(S) \not\subseteq R$.  Then, since $S$ is a minimal
extension of $R$, the subring $R +
\NST(S)$ of $S$ must equal $S$.  For $s \in S$, we can
find $r \in R$, $t \in \NST(S)$ such that $s = r+t$, and the
image of $r$ in $R/(\NST(S) \cap R)$ is uniquely determined by $s$.  It is
straightforward to see that the map sending
$s$ to the image of $r$ is a surjective $R$-homomorphism 
$S \rightarrow R/(\NST(S) \cap R)$ with kernel $\NST(S)$.  In
particular, $R/(\NST(S) \cap R)$ is $R$-isomorphic to $S/\NST(S)$.  
Since $S/\NST(S)$ is semiprime, we conclude that $\NST(S) \cap R
\supseteq \NST(R)$.  The reverse containment holds in general, by
\cite[Exercise 10.18(a)]{fc}, and so $\NST(R) = \NST(S) \cap R$.     
\end{proof}

Returning to the main subject of this section, 
we will study the ideal theory of ideal extensions.  
In particular, we will use the ideal theory 
to obtain information about the following cancellation problem: 
does the $R$-isomorphism class of $E(R,I)$ determine the $R$-isomorphism class of $I$?  
Fundamentally, this is a question about the ideals of $E(R,I)$, specifically regarding the ideals $I'$ of $E(R,I)$ 
which intersect $R$ trivially, and for which $R + I' = E(R,I)$.  When $I$ is a minimal $R$-rng, we will show 
that the above question can be answered in the affirmative. 

The following lemma characterizes the relevant ideals of $E(R,I)$, 
relating them to the set $\Hom_R(I,R)$.      

\begin{lemma} \label{posers}
Let $R$ be a ring, and let $I$ be an $R$-rng.  
Given $\varphi \in \Hom_R(I,R)$, define $I_{\varphi} = \{(\varphi(i),-i): i \in I\}$.  
For each $\varphi \in \Hom_R(I,R)$, $I_{\varphi}$ is an ideal of
$E(R,I)$ for which $R \oplus I_{\varphi} = E(R,I)$ as $($an internal
direct sum of$)$ abelian groups.  Conversely, if $I'$ is an ideal of $E(R,I)$ for which  $R \oplus I' = E(R,I)$ as abelian groups, then there exists a unique map $\varphi \in \Hom_R(I,R)$ such that $I' = I_{\varphi}$.    
\end{lemma}
\begin{proof}
Given $\varphi \in \Hom_R(I,R)$, it is easy to verify that
$I_{\varphi}$ is an ideal of $E(R,I)$, and we will only outline the
argument, leaving the details (which are similar to those appearing
two paragraphs below) to the reader.  The fact
that $\varphi$ is an $(R,R)$-bimodule homomorphism shows that
$I_{\varphi}$ is an additive subgroup of $E(R,I)$, and that
$I_{\varphi}$ is preserved by multiplication on either side by $R$.  
It remains only to show that $I_{\varphi}$ is stable under
multiplication by $I$ or by $I_{\varphi}$, 
since $E(R,I) = R + I = R + I_{\varphi}$.  
Computations similar to those found two paragraphs below can be used to
show either one of these statements.    

For the converse, let $I'$ be an ideal of $E(R,I)$ for which $R \oplus I' = E(R,I)$ as
abelian groups.  Since $R+I' = E(R,I)$, 
for each $i \in I$, there must be some $r \in R$ such that $(r,i) \in
I'$.  Since $I' \cap R = 0$, there is in fact a unique such $r$, since
$(r,i),(r',i) \in I'$ implies that $(r-r',0) \in R \cap I' = 0$.  Thus,
sending $i \in I$ to the unique $r \in R$ for which $(r,-i) \in I$
defines a map $\varphi: I \rightarrow R$, and 
$I' = \{(\varphi(i), -i): i \in I\}$.  

We claim that $\varphi: I \rightarrow R$ is an $R$-homomorphism
(i.e., a homomorphism of $R$-rngs).  First, let us show that $\varphi$
respects the $(R,R)$-bimodule structure.  
Suppose that $i,j \in I$ and let $r \in R$.  
Then, 
$(\varphi(i),-i) + (\varphi(j),-j) = (\varphi(i+j),-(i+j))$ since the
left-hand side is an element of $I'$ which has second component $-(i+j)$,
and the right-hand side is, by definition, the unique such element.  
We conclude that $\varphi(i) + \varphi(j) = \varphi(i+j)$.  
Similarly, $r(\varphi(i),-i) = (\varphi(ri),-ri)$ and
$(\varphi(i),-i)r = (\varphi(ir),-ir)$, from which we conclude that $\varphi$
is an $(R,R)$-bimodule homomorphism.  
Finally, observe that $(\varphi(i), -i)(\varphi(j),-j) =
(\varphi(i)\varphi(j), -i \varphi(j) - \varphi(i) j + ij)$, from
which we conclude that $\varphi( i \varphi(j) + \varphi(i)j - ij) =
\varphi(i) \varphi(j)$.  Using the fact that $\varphi$ is an
$(R,R)$-bimodule homomorphism, and the fact that
$\varphi(i),\varphi(j) \in R$, the left-hand side is
$\varphi(i)\varphi(j) + \varphi(i)\varphi(j) - \varphi(ij)$;
comparing with the right-hand side, we conclude that $\varphi(ij) =
\varphi(i)\varphi(j)$.  Therefore,  $\varphi \in \Hom_R(I,R)$.
\end{proof}

Using Lemma~\ref{posers}, we can quickly describe a condition on $\varphi$ under which
$I_{\varphi}$ and $I$ must be $R$-isomorphic, and from this we will
obtain cancellation, in the sense described above, for minimal $R$-rngs.      
\begin{lemma}  \label{suffiso}
  Let $R$ be a ring, and let $\varphi: I \rightarrow R$ be a homomorphism of $R$-rngs for
  which $\varphi(i)j = ij = i \varphi(j)$ for all
  $i,j \in I$.  Then the map $\Phi: I \rightarrow I_{\varphi}$, defined
  by $\Phi(i) = (\varphi(i),-i)$ is an $R$-isomorphism.  
  In particular, the above holds whenever $\varphi: I \rightarrow R$
  is an injective $R$-homomorphism.  
\end{lemma}
\begin{proof}
  The map $\Phi$ is clearly a bijective $(R,R)$-bimodule homomorphism, so we
  need only show that $\Phi$ is an $R$-rng homomorphism.  
  Thus, let $i,j \in I$, and note that  
  $\Phi(i)\Phi(j) = (\varphi(i),-i)(\varphi(j),-j) = (\varphi(ij), -i
  \varphi(j) - \varphi(i)j + ij)$.  By assumption, the right-hand side
  is $(\varphi(ij), -ij) = \Phi(ij)$, as desired.  

  For the last statement, note that if $\varphi$ is injective, then 
  $\varphi(i) j, ij$, and $i\varphi(j)$ must all agree, since
  $\varphi$ sends each to $\varphi(i)\varphi(j)$.  
\end{proof}

\begin{prop}  \label{idealcancel} Let $R$ be a ring, and let $I$ be an $R$-rng for which the only noninjective
  $R$-homomorphism $I \rightarrow R$ is the zero map $($in
  particular, this holds when $I$ is a minimal $R$-rng$)$.  Then, the
  $R$-isomorphism class of $E(R,I)$ determines the $R$-isomorphism
  class of $I$.  
\end{prop}
\begin{proof}
  Suppose that $I$ and $I'$ are $R$-rngs for which 
  $E(R,I)$ and $E(R,I')$ are $R$-isomorphic.  Under such
  an isomorphism $I'$ is $R$-isomorphic to some ideal $L$ (which we
  may view as an $R$-rng) of $E(R,I)$
  for which $R \oplus L = E(R,I)$ as abelian groups.  
  By Lemma~\ref{posers}, we conclude that $L = I_{\varphi}$ for some
  $R$-homomorphism $\varphi: I \rightarrow R$.  By hypothesis,
  either $\varphi$ is zero, in which case $L = I$; or else $\varphi$
  is injective, and $L = I_{\varphi}$ is $R$-isomorphic to $I$ by Lemma~\ref{suffiso}.  
  In any case, we conclude that $I'$ is $R$-isomorphic to $I$.  
\end{proof}

\section{The ideal theory of $E(R,I)$}\label{idealtheory}

In this section we will give a full description of the ideals of
$E(R,I)$, for an arbitrary minimal $R$-rng $I$ satisfying $I^2 \ne 0$.  We will also determine
when $E(R,I)$ is (semi)prime, and finally, we will discuss three
mutually exclusive classes of ideal extensions, which stratify the
ideal-theoretic behavior of an ideal extension (in particular,
controlling whether it is (semi)prime).  

Ideas similar to those used in the proof of Lemma~\ref{posers} 
can be used to give a full description of the ideals of $E(R,I)$, for an arbitrary $R$-rng $I$ satisfying $I^2 \ne 0$.  
Since we are concerned primarily with minimal extensions here, we will only include the 
description in that case, leaving the general case (which is only
slightly more technical) as an exercise to the interested reader.   
\begin{prop}  \label{idealdescription} Let $R$ be a ring, and let $I$ be a minimal $R$-rng with $I^2 \ne 0$.  The following is a complete list of the ideals of $E(R,I)$.  
\begin{enumerate}
	\item  $A \oplus 0$, where $A$ is an ideal of $R$ contained in $\, \ann_R(I)$. 
	\item  $A \oplus I$, where $A$ is an ideal of $R$.
	\item  $\{(a, -i): i \in I, a \in R, \text{ such that } a+Z = \varphi(i)\}$, where $Z
	\subseteq \ann_R(I)$ is an ideal of $R$, and $\varphi: I
	\rightarrow R/Z$ is a nonzero $R$-homomorphism.  
\end{enumerate}
The first type consists of all those ideals contained in $R \oplus 0$;
the second consists of all ideals which contain $I$, and the third
type consists of all other ideals $($those which neither contain $I$ nor
are contained in $R \oplus 0$$)$.  The last collection of ideals is
nonempty
if and only if $\Hom_R(I,R/\ann_R(I)) \ne 0$.
\end{prop}
\begin{proof}
If $A$ is an ideal of $R$ contained in $\ann_R(I)$, it is straightforward to see that $A \oplus 0$ is an ideal of $E(R,I)$.
Conversely, if $I'$ is an ideal of $E(R,I)$ with $I' \subseteq R$, then $I'I$ and $II'$ are both contained in $I \cap I' \subseteq I \cap R = 0$.  It follows that $I' \subseteq \ann_R(I)$, so $I' = A \oplus 0$, where $A$ is an ideal of $R$ contained in $\ann_R(I)$.    

Next, if $A$ is any ideal of $R$, then $A \oplus I$ is an ideal of $E(R,I)$ which contains $I$.  
Conversely, suppose that $I'$ is an ideal of $E(R,I)$ which contains $I$.  Let $A$ be the set of all 
$r \in R$ for which $(r,i) \in I'$ for some $i \in I$.  Clearly $A$ is an ideal of $R$, and $I' = A \oplus I$, as claimed.    

Finally, suppose that $K$ is an ideal of $E(R,I)$ which does not contain
$I$ and is not contained in $R$.
Consider the set $C$ of all second coordinates of elements of $K$; that is, the set 
of $i \in I$ for which $(r,i) \in K$ for some $r \in R$.  It is
clear that $C$ is an $(R,R)$-subbimodule of $I$, and that $C$ is
nonzero, since $K$ is not contained in $R$.  
Now suppose that $(r,i),(r',i') \in K$.  The equation
$$(r,i)(r',i') - (r,i)r' - r(r',i') = (-rr', ii')$$
establishes that $C$ is an $R$-subrng of $I$, so $C = I$ by
minimality of $I$.  
Now, let $Z = \{r \in R: (r,0) \in K\}$, which is an ideal of $R$.  
Define a map $\varphi: I \rightarrow R/Z$ as follows.  Given $i \in
I$, we may find some $r \in R$ such that $(r,-i) \in K$.  Define
$\varphi(i)$ to be the image of such an $r$ in $R/Z$.  To see that
$\varphi$ is well-defined, it suffices to note that if $(r,-i),(r',-i)
\in K$, then $(r-r',0) \in K$, and hence $r-r' \in Z$.  To see that
$\varphi$ is an $R$-homomorphism, it suffices to note that if
$(r,-i),(r',-i') \in K$, then $(r+r',-(i+i')),(rr',-ri'), (rr',-ir')$,
and $(rr',-ii')$ are in $K$ (the last established by the equality
above); the membership of the first three in $K$ shows (reducing modulo $Z$) that $\varphi$ is an $(R,R)$-bimodule
homomorphism, while that of the last shows (reducing modulo $Z$) that $\varphi(ii') =
\varphi(i)\varphi(i')$.  Thus $\varphi$ is an $R$-homomorphism, and it is now easy to check that $K = \{(a,-i): i \in I, a \in R,
a + Z = \varphi(i)\}$.  Indeed, for each $i \in I$, we can find some
$a \in R$ such that $(a,-i) \in K$.  But $Z \oplus 0 \subseteq K$,
and hence $(a+Z,-i) \subseteq K$.  It follows that $a+Z = \varphi(i)$.  
Since $K$ does not contain $I$, $\varphi \ne 0$, since otherwise, $K =
Z \oplus I$, which does contain $I$.  Finally, $I$ is $R$-isomorphic
to its image in $R/Z$ under $\varphi$ (which is nonzero),
and hence $Z \subseteq \ann_R(I)$, since $Z$ annihilates $R/Z$.    

Conversely, let $Z \subseteq \ann_R(I)$ be an ideal of $R$ and let
$\varphi: I \rightarrow R/Z$ be a nonzero $R$-homomorphism.  Consider
the subset $K = \{(a,-i): i \in I, a \in R, a+Z = \varphi(i)\}$ of
$E(R,I)$.  Since $\varphi$ is an $R$-homomorphism, it is
straightforward to see that $K$ is an $(R,R)$-subbimodule of $E(R,I)$.
Since $R+I = E(R,I)$, to finish showing that $K$ is an ideal of $E(R,I)$, it will suffice to
show that $IK$ and $KI$ are contained in $K$.  Thus, let $i' \in I$ and 
$(a,-i) \in K$.  Then, $(a,-i)(0,i') = (0,ai' -ii')$.  Since
$\varphi$ is an $R$-homomorphism, we have 
$\varphi(ai' - ii') = a \varphi(i') - \varphi(i)\varphi(i')$.  But, 
by assumption $\varphi(i) = a + Z$, so $\varphi(ai' - ii') = 0 + Z$,
which is to say that $(0,ai' - ii') \in K$.  It follows that $KI
\subseteq K$, and similarly $IK \subseteq K$, so we conclude that $K$
is an ideal of $E(R,I)$.  
Finally, it is clear that $K$ is not contained in $R$, and the fact that
$\varphi$ is nonzero ensures that $K$ does not contain $0 \oplus I$.    

For the final statement, suppose that there is an ideal of the form (3).  	
Then, there is an ideal $Z \subseteq \ann_R(I)$ of $R$ such that
$I$ is $R$-isomorphic to a minimal ideal of $R/Z$.  Under the further
quotient map $R/Z \rightarrow R/\ann_R(I)$, the image of $I$ in $R/Z$ cannot be sent to
zero, since $I^2 \ne 0$.  Thus, the $R$-isomorphic image of $I$ in $R/Z$
is not contained in $\ann_R(I)$.  It follows that
$\Hom_R(I,R/\ann_R(I)) \ne 0$.  Conversely, if there is a nonzero element in $\Hom_R(I,R/\ann_R(I))$, we may use such a map with $Z = \ann_R(I)$ to produce an ideal of type (3).    
\end{proof}   


We next characterize when $E(R,I)$ is (semi)prime.

\begin{prop}  \label{semiprimeoversemiprime}  Let $R$ be a ring, and let $I$ be a minimal $R$-rng.  
Then, $E(R,I)$ is semiprime if and only if $R$ is semiprime and $I^2
\ne 0$.  
\end{prop}
\begin{proof} 
For the forward implication, suppose that $E(R,I)$ is semiprime.  Since $I$
is a nonzero ideal of $E(R,I)$, we must have $I^2 \ne 0$, and hence
$\ann_R(I)$ is prime, by Lemma~\ref{minimalann}.  Let $N$ be a
nilpotent ideal of $R$.  Then clearly $N \subseteq \ann_R(I)$, since
the image of $N$ in the prime ring $R/\ann_R(I)$ is nilpotent, hence
zero.  If $K \subseteq \ann_R(I)$ is an ideal of $R$, then $K \oplus
0$ is an ideal of $E(R,I)$, and hence we conclude that $N \oplus 0$ is
an ideal of $E(R,I)$.  But $N \oplus 0$ is nilpotent, so we conclude
that $N = 0$.  It follows that $R$ is semiprime.  

For the reverse implication, suppose that $E(R,I)$ is not semiprime and that  
$A$ is a nonzero ideal of $E(R,I)$ with $A^2 = 0$.  
Since $R$ is semiprime and $(R \cap A)^2 = 0$, we have $R \cap A = 0$.  But $A$ is nonzero, so we must have 
$R \oplus A = E(R,I)$ as abelian groups, since $E(R,I)$ is a minimal extension of $R$, by Lemma~\ref{produce}.  
We conclude that $E(R,I)$ and $E(R,A)$ are $R$-isomorphic, and hence 
$I$ and $A$ are $R$-isomorphic by Proposition~\ref{idealcancel}.  Since $A^2 = 0$, we must have $I^2 = 0$ as well.  
\end{proof}

The following is an analogous statement for prime rings. 
\begin{prop}  \label{primeidealext} Let $R$ be a ring, let $I$ be a
  minimal $R$-rng, and let $E = E(R,I)$ be the associated $($minimal$)$ ideal
  extension.  Then, the following are equivalent.
\begin{enumerate}
	\item  $E$ is prime $($and is subdirectly irreducible$)$. 
	\item  $\ann(I_{E}) = \ann({}_{E}I)= 0$.
	\item  $I^2 \ne 0$, $\, \ann_R(I) = 0$, and $\, \Hom_R(I,R) = 0$.  
\end{enumerate}
In particular, if $E(R,I)$ is prime, then $R$ must be prime.  
\end{prop}
\begin{proof}
The implication $(1) \Rightarrow (2)$ is clear, since $I \ne 0$.  
To prove $(2) \Rightarrow (3)$, suppose that $\ann(I_E) = \ann({}_E I)
= 0$.  
This implies, in particular, that $I^2 \ne 0$, and $\ann_R(I) = 0$.  We
claim that if $0 \ne \varphi \in \Hom_R(I,R)$, then $I_{\varphi} I = 0$ (see Lemma~\ref{posers} for the definition of $I_{\varphi}$).  Indeed, if $(\varphi(i),-i) \in I_{\varphi}$, and $(0,i') \in I$, we have 
$(\varphi(i),-i)(0,i') = (0, \varphi(i)i' - ii')$.  Since $I$ is minimal, $\varphi$ must be injective, and hence 
$\varphi(i)i' = ii'$ (as in Lemma~\ref{suffiso});  we conclude that 
$I_{\varphi} I = 0$, which contradicts $\ann({}_E I) = 0$.  It follows that $\Hom_R(I,R) = 0$.  

To prove $(3) \Rightarrow (1)$ and the final claim, suppose that $I^2
\ne 0$, $\ann_R(I) = 0$, and $\Hom(I,R) = 0$.  Then $R$ is prime,
since $\ann_R(I) = 0$ is a prime ideal of $R$, by Lemma~\ref{minimalann}.  Suppose that $A$ and $B$ are nonzero ideals of $E(R,I)$ for which $AB = 0$.  
Since $(R \cap A)(R \cap B) = 0$ and $R$ is prime, we conclude that
either $R \cap A = 0$ or $R \cap B = 0$.  Without loss of generality,
$R \cap A = 0$.  By minimality of $E(R,I)$, $R+A = E(R,I)$, and from
Lemma~\ref{posers}, we conclude that $A = I_{\varphi}$ for some
$\varphi \in \Hom_R(I,R)$.  But $\Hom_R(I,R) = 0$, so $\varphi = 0$, and hence $A = I$.  Now, $(B \cap R) \subseteq \ann(A_R) = \ann(I_R) = \ann_R(I) = 0$.  We conclude that $B \cap R = 0$, and since $B$ is nonzero, we conclude that $B = I$, as we did with $A$.  
But now, $AB = I^2$ is nonzero, a contradiction.  We conclude that $E(R,I)$ is prime, as desired.   
\end{proof}

\begin{cor}  \label{annideals}
Let $R$ be a ring, and let $I$ be a minimal $R$-rng.  
\begin{enumerate}
	\item $\ann_R(I) \oplus 0$ is a semiprime ideal of $E(R,I)$ if and only if $I^2 \ne 0$.  
	\item $\ann_R(I) \oplus 0$ is a prime ideal of $E(R,I)$ if and
	only if $I^2 \ne 0$ and $\, \Hom_R(I,R/\ann_R(I)) = 0$.  In this case, 
	$E(R,I) / (\ann_R(I) \oplus 0) \cong E(R/\ann_R(I),I)$ is a subdirectly irreducible prime ring.    
\end{enumerate}
\end{cor}
\begin{proof}
We have the $R$-ring homomorphism $E(R,I)/(\ann_R(I) \oplus 0) \cong E(R/\ann_R(I),I)$, where $I$ is viewed as an 
$(R/\ann_R(I))$-rng.  We also note that $R/\ann_R(I)$ is semiprime by Lemma~\ref{minimalann}.  The first
statement now follows 
from Lemma~\ref{semiprimeoversemiprime}, and the second follows from Lemma~\ref{primeidealext}, once we observe that $\ann_{R/\ann_R(I)}(I) = 0$.          
\end{proof}

\begin{remark} \label{setupthreetypes} In light of Propositions~\ref{idealdescription},
  \ref{semiprimeoversemiprime} and \ref{primeidealext}, it is natural to group minimal ideal
  extensions $E(R,I)$ into three types, based on the following
  properties of $I$ which control the ideal-theoretic behavior of $E(R,I)$.       
\begin{enumerate}
  \item  $I^2 = 0$,
  \item  $I^2 \ne 0$ and $\Hom_R(I,R/\ann_R(I)) = 0$, 
  \item  $\Hom_R(I,R/\ann_R(I)) \ne 0$ (which forces $I^2 \ne 0$,
  since $R/\ann_R(I)$ is always semiprime; this in turn implies that
  $R/\ann_R(I)$ is prime by Lemma~\ref{minimalann}).  
\end{enumerate}
By Proposition~\ref{idealcancel}, $R$-rngs $I$ falling under different cases above produce non-$R$-isomorphic ideal extensions $E(R,I)$.  
\end{remark}

Recall that a ring $R$ is said to be subdirectly irreducible 
if it has a least nonzero ideal, called the {\it little ideal} (see \cite[Section 12]{fc}).  We will say that an ideal $P$ of $R$ is subdirectly irreducible if $R/P$ is a subdirectly irreducible ring; clearly, any maximal ideal is subdirectly irreducible.
Let us give a better description of the third type of minimal
$R$-rngs in the above list, those for which $\Hom_R(I,R/\ann_R(I)) \ne 0$.          
\begin{lemma} \label{subdirectprime}
Let $R$ be a ring.  The following are equivalent.  
\begin{enumerate}
\item  $I$ is a minimal $R$-rng with $\, \Hom_R(I,R/\ann_R(I)) \ne 0$.
\item  $R/\ann_R(I)$ is a subdirectly irreducible prime ring, and $I$ is $R$-isomorphic to its little ideal $($as an $R$-rng$)$.  
\end{enumerate}
\end{lemma}
\begin{proof}
Observe that a prime ring with a minimal nonzero ideal must be subdirectly irreducible.  
To see this, let $S$ be a prime ring with a minimal nonzero ideal $K$, and let $K'$ be any nonzero ideal of $S$.  
Then, $KK'$ must be nonzero, which implies that $K \cap K' \ne 0$.  By minimality of $K$, we must have $K \cap K' = K$, and 
hence $K \supseteq K'$.  Thus, every nonzero ideal of $S$ contains $K$, so $S$ is subdirectly irreducible.  

Now, suppose that $I$ is a minimal $R$-rng and that $\Hom_R(I,R/\ann_R(I)) \ne 0$.  
Then, $\ann_R(I)$ is a prime ideal of $R$, and the image of $I$ in $R/\ann_R(I)$ under any nonzero $R$-homomorphism $I \rightarrow R/\ann_R(I)$ is a minimal ideal of the prime ring 
$R/\ann_R(I)$.  It follows that $R/\ann_R(I)$ is a subdirectly irreducible prime ring, and $I$ is $R$-isomorphic to its little ideal (viewed as an $R$-rng).  
The converse is clear.  
\end{proof}
Thus, the $R$-isomorphism classes of minimal $R$-rngs $I$ for which $\Hom_R(I,R/\ann_R(I)) \ne 0$
are in one-to-one correspondence with subdirectly irreducible prime ideals of $R$.  

The following lemma characterizes each of the three types of minimal $R$-rngs $I$
based on properties of the annihilator of $I$ in $E(R,I)$.  
\begin{lemma} \label{thethreetypes}
Let $R$ be a ring, and let $I$ be a minimal $R$-rng, so that $E(R,I)$ is a minimal extension of $R$.  Then the following hold.  
\begin{enumerate}
   \item  $I^2 = 0$ if and only if $\, \ann_{E(R,I)}(I) \supseteq I$.  
   \item  $I^2 \ne 0$ and $\Hom_R(I,R/\ann_R(I)) = 0$, if and only if $\, \ann_{E(R,I)}(I) \subseteq R$.   
   \item  $\Hom_R(I,R/\ann_R(I)) \ne 0$ if and only if
   $\, \ann_{E(R,I)}(I)$ intersects $I$ trivially and is not contained in
   $R$.  
\end{enumerate}
\end{lemma}
\begin{proof}
The first statement is clear.  For each of the other two statements, 
we may assume that $\ann_R(I) = 0$, by
passing to the quotient $E(R,I)/(\ann_R(I) \oplus 0)$, which is
$R$-isomorphic to $E(R/\ann_R(I),I)$.  Under the quotient map
$\ann_{E(R,I)}(I)$ is sent to $\ann_{E(R/\ann_R(I),I)}(I)$ (since the
ideal $I$ intersects the kernel of the map trivially); in particular,
the statements regarding $E(R,I)$ and containment in $R$ are
equivalent to the corresponding statements with ``equals 0'' replacing
containment in $R$.  

Let us prove statement (2).  For the forward implication,
suppose that $\Hom_R(I, R) = 0$, $I^2 \ne 0$, and $\ann_R(I) = 0$.  
If $A = \ann_{E(R,I)}(I) \ne 0$, then by Lemma~\ref{posers} we must
have $A = I_{\varphi}$ for some $\varphi \in \Hom_R(I,R)$, since $A \cap R = \ann_R(I) = 0$.    
But $\Hom_R(I,R) = 0$, so $\varphi = 0$, and hence $A = I$, which implies that $I^2 = 0$, which is a contradiction. 
We conclude that $\ann_{E(R,I)}(I) = 0$, as desired.  Conversely, suppose that
$\ann_{E(R,I)}(I) = 0$.  Then $I^2 \ne 0$, and we claim that $\Hom_R(I,R) = 0$.
Indeed, supoose that $0 \ne \varphi \in \Hom_R(I,R)$.  Then $I_{\varphi} I = 0 = I I_{\varphi}$, as in the proof 
of Proposition~\ref{primeidealext}, which contradicts the fact that $\ann_{E(R,I)}(I) = 0$.  We conclude that 
$\Hom_R(I,R)= 0$ and $I^2 \ne 0$, as desired.  

The third statement follows from the second, since $\Hom_R(I,R/\ann_R(I)) \ne 0$ implies that $I^2 \ne 0$, which, in turn, implies that $\ann_{E(R,I)}(I) \cap I = 0$.
\end{proof}

\section{Central extensions} \label{centralsection}

In this section, we will examine a class of minimal ideal extensions that
behave fundamentally like commutative minimal extensions.  We start
with a definition.        
\begin{defn}  
Let $R$ be a subring of a ring $S$.  We say that $S$ is a {\emph{central $($ring$)$
extension}} of $R$ if as a (left) $R$-module, $S$ is generated by
elements of $\, \C_S(R)$, the centralizer of $R$ in $S$.  
\end{defn}  

We note that if $S$ is a central minimal extension of a ring $R$, and  
$s \in \C_S(R) \setminus R$ (which is nonempty by hypothesis), then $S
= R[s]$ (the subring of $S$ generated by $R \cup \{s\}$).  
In particular, a central minimal extension of a commutative ring must
be commutative.  
In general, central extensions of commutative rings need not be commutative (e.g., the rational quaternions 
over $\Q$).  Also note that if $S$ is a central minimal extension of
$R$, then $R$ need not belong to $Z(S)$, the center of $S$ (for instance,
if $R$ is not commutative).       

\begin{lemma}  \label{centralstuff}  Let $R$ be a ring, and let $I$ be a minimal $R$-rng.  
Then $E(R,I)$ is a central extension of $R$ if and only if $\, \C_I(R) = \{ x \in
I : x r = r x \text{ for all } r \in R\}$ is nonzero.  Moreover,
$\, \C_I(R) \subseteq Z(E(R,I))$.       
\end{lemma}
\begin{proof}
For the forward implication,
we are given that $(r,i) \in E(R,I)$ centralizes $R$, for some $r \in
R$, and $0 \ne i \in I$.  Then, for any $r' \in R$, we have $(r',0)(r,i) =
(r,i)(r',0)$, from which we conclude that $r'i = ir'$, and hence $0
\ne i \in \C_I(R)$.   

The converse and the final claim follow from the fact that if $0 \ne i \in \C_I(R)$, then $E(R,I) = R[i]$.
\end{proof}

To characterize when $E(R,I)$ is a central extension, we will need 
a lemma, which uses the ideas behind the proof of Brauer's Lemma
(e.g., cf.\ \cite[Lemma 10.22]{fc}).  
\begin{lemma}  \label{Brauer} Let $R$ be a ring, and let $I$ be a
  minimal $R$-rng with $I^2 \ne 0$.   
  Then, the following are equivalent.  
  \begin{enumerate}
    \item  $I \cap Z(R) \ne 0$.
    \item  $I$ contains a nonzero central idempotent of $R$. 
    \item  $I$ is a direct summand of $R$ $($as $2$-sided ideals$)$.  

  \end{enumerate}
\end{lemma}  
\begin{proof}
The equivalence of the second and third statements is standard, and
the implication $(2) \Rightarrow (1)$ is obvious.  
For the implication $(1) \Rightarrow (2)$, suppose that $0 \ne x
\in I \cap Z(R)$.
We must have $Rx = I$, since $Rx$ is a nonzero $2$-sided ideal of $R$
contained in the minimal ideal $I$.  Note that $Rx^2 = I^2 \ne
0$, from which we conclude that $x^2 \ne 0$.  Now, $Ix$ is a $2$-sided
ideal of $R$ contained in $I$, and is nonzero since it contains $x^2$,
so we conclude that $Ix = I$.  Thus, we may find $e \in I$ such that
$ex = x$ (necessarily, $e \ne 0$).  We
note that $A = \{i \in I : ix = 0\}$ is a proper sub-$R$-rng of $I$ (since
$e \not\in A$), and hence $A = 0$, by minimality.  Therefore, $e(e-1)x = 0$ implies
that $e(e-1) = 0$, from which we conclude that $e^2 = e$.  Moreover,
if $r \in R$, then $(re - er)x = 0$.  But $re - er \in I$, so we
conclude that $re - er \in A = 0$.  Thus, $e$ is a
nonzero central idempotent of $R$ contained in $I$.  
\end{proof}

\begin{cor}  \label{primecenter} Let $R$ be a subdirectly irreducible prime ring with
  little ideal $I$.  Then, $I \cap Z(R) \ne 0$ if and only if $R$ is
  simple.  
\end{cor} 
\begin{proof}
For the forward implication, we know from
Lemma~\ref{Brauer} that $I$ contains a nonzero central idempotent $e$
of $R$.  If $e \ne 1$, then $R$ is a nontrivial direct product of two
rings, so $R$ is not prime.  We conclude that $1 = e \in I$.  Therefore, the
little ideal of $R$ is the whole ring, and hence $R$ is simple.    
If $R$ is simple, then $1 \in R = I$ which establishes the reverse
implication.  
\end{proof}

The next result gives several descriptions of when $E(R,I)$ is a central extension of $R$, in the case where $I$ is a minimal $R$-rng satisfying $I^2 \ne 0$.

\begin{prop} \label{maincentral}
Let $R$ be a ring, and let $I$ be a minimal $R$-rng.  The following are
equivalent. 
\begin{enumerate}
\item $I$ is $R$-isomorphic to $R/\ann_R(I)$.  
\item $I$ is a ring.
\item $I$ has a nonzero central idempotent.  
\item $\Hom_R(R,I) \ne 0$.  
\item $\Hom_R(I,R/\ann_R(I)) \ne 0$ and $\, \ann_R(I)$ is a maximal ideal
  of $R$.  
\item $E(R,I)$ is a central extension of $R$ and $I^2 \ne 0$.  
\end{enumerate}
\end{prop}
\begin{proof}
The implications $(1) \Rightarrow (2) \Rightarrow (3)$ are
trivial, even without assuming that $I$ is minimal.  

The implication $(3) \Rightarrow (4)$ also does not require
minimality.  Let $e \in I$ be a nonzero central idempotent (note that
we are not assuming that $e$ commutes with the action of $R$).  
Consider the map $f: R \rightarrow I$ defined by $f(r) = re$.  
It is clear that $f$ is a left $R$-module homomorphism, but we
claim that $f$ is actually an $R$-homomorphism.  
Note that (using only the fact that $e$ commutes with elements of
$I$) for all $r,s \in R$,  $$f(rs) =
(rs)e = (rs)ee = r(se)e = re(se) = r(es)e = re(es) = res.$$   
These equations imply that 
$f(rs) = f(r)f(s)$ and $f(rs)
= f(r)s$, from which it follows that $f$ is an $R$-homomorphism.    

To prove the implication $(4) \Rightarrow (1)$, 
let $0 \ne f \in \Hom_R(R,I)$, and set $e =
    f(1)$, which is a central idempotent of $I$ that commutes with
    the action of $R$.  Clearly $f(R) = eR =
    Re$ is a nonzero $R$-subrng of $I$, and hence the minimality of
    $I$ implies that $f(R) = I$.  It follows that $I$ is $R$-isomorphic to
    $R/\ker(f)$.  Finally, if $x \in R$, then $f(x) = xf(1) = f(1) x$,
    so clearly $\ker(f) = \ann_R(e) = \ann_R(I)$.  

Statements $(1)$ and $(5)$ are equivalent, since $R/\ann_R(I)$ is a
minimal $R$-rng if and only if $\ann_R(I)$ is a maximal ideal of $R$, and any nonzero
$R$-homomorphism between minimal $R$-rngs is an $R$-isomorphism.  

Next, $(1) \Rightarrow (6)$, since $1 \in R/\ann_R(I) \cong I$
centralizes $R$, and clearly $I^2 \ne 0$.  Finally, we will show that $(6) \Rightarrow (3)$.  Suppose that $x
\in E(R,I) \setminus R$ centralizes $R$.  By Lemma~\ref{centralstuff}, 
$0 \ne \C_I(R) \subseteq Z(E(R,I))$.  But $I$ is a minimal ideal of
$E(R,I)$ satisfying $I^2 \ne 0$, and $I \cap Z(E(R,I)) \ne 0$, so
Lemma~\ref{Brauer} 
implies that $I$ contains a central idempotent of
$E(R,I)$, as desired.  
\end{proof}

\begin{remark}  Of the three types of minimal $R$-rngs described in
  Remark~\ref{setupthreetypes}, by Proposition~\ref{maincentral}, 
  ideal extensions corresponding to those of the second type are never
  central, and an extension of the third type is central precisely when $\ann_R(I)$
  is a maximal ideal $($rather than merely a subdirectly irreducible
  prime ideal $P$$)$.
\end{remark} 

We are now ready to prove a criterion for recognizing when an arbitrary minimal ideal extension is central.

\begin{prop}  \label{centralchar}  Let $R$ be a ring, and let $I$ be a
  minimal $R$-rng.  Then $E(R,I)$ is a central extension of $R$ if and
  only if $I$ is isomorphic to $R/\ann_R(I)$ as an
  $(R,R)$-bimodule.  In particular, this implies that $\ann_R(I)$ is a
  maximal ideal.  
\end{prop}  
\begin{proof}  
For the forward implication, suppose that $E(R,I)$ is a central
extension of $R$.  If $I^2 \ne 0$, then the conclusion follows from
Proposition~\ref{maincentral}.  Thus, suppose that $I^2 = 0$.  By
Lemma~\ref{centralstuff}, we can find a nonzero $i \in I \cap Z(E(R,I))$.  The map $f: R
\rightarrow I$, defined by $f(r) = ri$, is an $(R,R)$-bimodule
homomorphism, and $f(R)$ is a nonzero $(R,R)$-subbimodule of $I$.  
Since $I^2 = 0$, $f(R)$ is a nonzero $R$-subrng of $I$, and hence, by
minimality, $f(R) = I$.  Since $\ker(f) = \ann_R(i) = \ann_R(I)$, it
follows that $I$ is isomorphic to $R/\ann_R(I)$ as an
$(R,R)$-bimodule.  
The reverse implication is clear from Lemma~\ref{centralstuff}, 
since the image of $1$ in $R/\ann_R(I)$ (rather, its image in $I$)
centralizes $R$.  The final statement is immediate.
\end{proof}

We conclude this section with a few observations regarding the
behavior of the prime radical under central extensions.  
We observed earlier that non-semiprime rings can have semiprime
minimal ring extension (for instance, $\M_2(k)$ over $\T_2(k)$, for a
field $k$);  this phenomenon
does not persist in the central case.  
\begin{lemma}  \label{semiprimeovercentral} If $S$ is a central extension $($not necessarily minimal$)$
  of a ring $R$ and $S$ is
  $($semi$)$prime, then $R$ is $($semi$)$prime.  
\end{lemma} 
\begin{proof}  This is essentially the same proof as that of
  \cite[Exercise 10.18A(2)]{fc}.  Let $R$ be a ring, let $S$ be a
  prime central
  extension of $R$, and let $X$ be an $R$-centralizing set which generates
  $S$ as a left $R$-module.  Thus, $S = R \langle X \rangle$.       
  Suppose that $aRb = 0$ with $a,b \in R$.  
  Then, $aSb = a(R\langle X \rangle)b = aRb\langle X \rangle = 0$.  Since $S$ is prime, either $a=0$
  or $b=0$.  
Thus, $aRb = 0$ implies that $a=0$ or $b=0$, so $R$ is prime. 
The semiprime case is similar, and is left to the reader.  
\end{proof}
Lemma~\ref{semiprimeovercentral} implies that for a
minimal central extension $S$ of $R$, if $\NST(S) \subseteq R$, then $\NST(S)
= \NST(R)$.  To prove this, note that $S/\NST(S)$ is a central semiprime minimal
extension of $R/\NST(S)$, and apply 
Lemma~\ref{semiprimeovercentral}.
In particular, this, together with Lemma~\ref{partialreduction},
provides a dichotomy for central extensions analogous to that found in
\cite[Theorem 2.1]{D&S2}.

\section{Minimal extensions of prime rings} \label{primesection}

As we saw in Proposition~\ref{teaser}, every non-prime minimal
extension of a prime ring is an ideal extension.  In this section, we
will use Proposition~\ref{teaser} together with our results on ideal
extensions to classify the minimal extensions of arbitrary prime
rings.  Moreover, we will fit this together with the results of Dobbs and Shapiro on minimal commutative extensions of commutative domains.

We begin by recording the Dobbs-Shapiro classification of minimal
commutative extensions of a commutative domain 
(which follows the earlier Ferrand-Olivier classification of minimal
commutative extensions of fields, found in \cite{F&O}).  
\begin{thmnonum} \cite[Theorem 2.7]{D&S1} \label{dobbsshapiro} Let $R$ be a commutative domain.
  Up to $R$-isomorphism, 
every minimal commutative extension of $R$ is of exactly one of the
following forms.     
\begin{enumerate}
	\item[\tdom]  A domain that is a minimal extension of $R$.
	\item[\tred]  $R \times R/M$, for some $M \in \Max(R)$. 
	\item[\tnot]  $R \propto M$, for some $M \in \Max(R)$.  
\end{enumerate}
The maximal ideal $M$, where it appears, is determined by the
$R$-isomorphism type of the extension.  
\end{thmnonum}
The labels \tdom, \tred, and \tnot, refer to properties of the
associated extensions, namely, \ttdom refers to those extensions which
are domains, \ttred refers to those extensions which are reduced but
are not domains, and \ttnot refers to those extensions which are not
reduced.  In particular, from this it is clear that the type (\tdom,
\tred, or \tnot) of the extension is determined uniquely by the isomorphism type.  

Since each non-prime minimal extension of a prime ring is an ideal
extension, by Proposition~\ref{teaser}, type-\ttred and type-\ttnot minimal extensions above must be ideal extensions
(and this is easy to see directly).  
Using our work on ideal extensions, we will produce a classification for
arbitrary minimal extensions of prime rings (which looks substantially
different from the above classification), and we will
specialize it to the central case, where the classification looks
almost identical to \cite[Theorem 2.7]{D&S1}. 

We are now ready to prove our main result for prime rings.  
\begin{thm}  \label{primeext} Let $R$ be a prime ring.  Then, up to $R$-isomorphism, every minimal 
extension of $R$ must be of exactly one of the following five forms.
\begin{enumerate}
	\item[\ta]  A prime minimal extension of $R$, all of whose nonzero ideals intersect $R$ nontrivially. 
	\item[\tb]  $E(R,I)$ for some minimal $R$-rng $I$ such that
	$\, \Hom_R(I,R) = 0$, $I^2 \ne 0$, and $\, \ann_R(I) = 0$.  
	\item[\tc]  $E(R,I)$ for some minimal $R$-rng $I$ such that $\,\Hom_R(I,R/\ann_R(I)) = 0$, $I^2 \ne 0$, and $\, \ann_R(I) \ne 0$.
	\item[\td] $E(R,I)$, where $I$ is a minimal ideal of $R/P$
	for some prime ideal $P$ of $R$ $($note that this implies that
	$P$ is subdirectly irreducible; see Lemma~\ref{subdirectprime}$)$.  
	\item[\te] $R \propto M$ for some simple $(R,R)$-bimodule $M$.  
\end{enumerate}
Extensions of the forms \tta and \ttb are prime; those of the forms \ttc and
\ttd are semiprime, but not prime; and those of the form \tte are not
semiprime.  
In each case where
they occur, $I$, $M$, and $P$ are unique, up to $R$-isomorphism,
$(R,R)$-bimodule isomorphism, and equality, respectively.

The labels are intended to mean:  \tta = prime; \ttb = prime, ideal extension; \ttc = semiprime, reducible; \ttd = semiprime, subdirectly
irreducible; \tte = not semiprime.     
\end{thm}
\begin{proof}
By Lemma~\ref{produce}, an extension of $R$ of one of the above forms is minimal.  Further, an extension of type \tta clearly cannot be an ideal extension.  
On the other hand, an extension which has a nonzero ideal that 
intersects $R$ trivially must be an ideal extension of $R$.  In 
particular, as we saw in Proposition~\ref{teaser}, every non-prime 
minimal extension of $R$ must be an ideal extension.  Moreover, 
the types \tb-\tte include all possible minimal ideal extensions.  
Indeed, case \tte is the collection of ideal extensions $E(R,I)$ for 
which $I^2 = 0$; cases \ttb and \ttc are those extensions where 
$\Hom_R(I,R/\ann_R(I)) = 0$, and case \ttd consists of those ideal 
extensions $E(R,I)$ for which $\Hom_R(I,R/\ann_R(I)) \ne 0$, by Lemma~\ref{subdirectprime}.  
Uniqueness of the relevant data follows from Proposition~\ref{idealcancel}, 
moreover this shows that the case is determined by the $R$-isomorphism 
type of the ideal extension, since all data involved are determined by the
$R$-isomorphism type of $I$ or $M$ (note that for type-\ttd extensions, $P =
\ann_R(I)$).   

Finally, Proposition~\ref{primeidealext} shows that $E(R,I)$ is prime 
if and only if $I^2 \ne 0$, $\ann_R(I) = 0$, and $\Hom_R(I,R) = 0$, 
so extensions of forms \tta and \ttb are prime, and extensions of the
other types are not prime.  
By Proposition~\ref{semiprimeoversemiprime}, extensions of type \ttc or
\ttd are semiprime, since in those cases $I^2 \ne 0$ (and $R$ is
semiprime).  In addition, extensions of type \tte fail to be semiprime.    
\end{proof}

\begin{prop}  \label{centralprime} Let $R$ be a prime ring.  Adopting the same labeling and notation as in Theorem~\ref{primeext}, we have the following.
\begin{enumerate}
	\item  Type-\ttb and type-\ttc extensions of $R$ are never central. 
	\item  An extension of type \ttd is a central
	extension if and only if $\, \ann_R(I)$ is a maximal ideal of $R$, in which case the extension 
	is $R$-isomorphic to $R \times R/\ann_R(I)$.
	\item  An extension of type \tte is a central
	extension if and only if $\, \ann_R(M)$ is a maximal ideal of $R$, in which case the extension 
	is $R$-isomorphic to $R \propto R/\ann_R(M)$.
\end{enumerate}
\end{prop}  
\begin{proof}  
The first statement follows immediately from Proposition~\ref{maincentral}, since type-\ttb and type-\ttc ideal extensions are of the form $E(R,I)$, where $\Hom_R(I,R/\ann_R(I)) = 0$ and $I^2 \ne 0$.  

By Proposition~\ref{maincentral}, a
type-\ttd central extension must be of the form $E(R,I)$, where $I$ is
$R$-isomorphic to $R/\ann_R(I)$ and $\ann_R(I)$ is a maximal ideal.  
The map sending $(r, \overline{s}) \in R \times R/\ann_R(I)$ to $(r, \overline{s-r}) \in E(R,R/\ann_R(I))$ is easily seen to be an $R$-isomorphism from $R \times R/\ann_R(I)$ (with $R$ embedded diagonally) to $E(R,R/\ann_R(I))$.  

Finally, type-\tte central minimal extensions are of the form $E(R,M)$, for an $R$-rng $M$ satisfying $M^2 =0$, such that $M$ is isomorphic to 
$R/\ann_R(M)$ as an $(R,R)$-bimodule, by Proposition~\ref{centralchar} (where $\ann_R(I)$ is a maximal ideal).  
But $E(R,M)$ is then $R$-isomorphic to $R \propto R/\ann_R(M)$.    
\end{proof}

The following corollary classifies the central minimal extensions of a prime ring, in a form similar  to that of Theorem~\ref{primeext} .  The special case where $R$ is commutative is precisely \cite[Theorem 2.7]{D&S1}.    

\begin{cor}  \label{succinctcentral}
Let $R$ be a prime ring.  Then, up to $R$-isomorphism, every central
minimal extension of $R$ must be of exactly one of the following three forms.
\begin{enumerate}
\item[\ta]  A prime minimal extension of $R$, all of whose nonzero
  ideals intersect $R$ nontrivially.  
\item[\td]  $R \times R/M$, where $M \in \Max(R)$.
\item[\te]  $R \propto R/M$, where $M \in \Max(R)$.
\end{enumerate}
The maximal ideal $M$, where it appears, is uniquely determined by the
$R$-isomorphism type of the extension.  
\end{cor}

\begin{proof}
The classification follows from Theorem~\ref{maincentral} and  Proposition~\ref{centralprime} .  The final statement follows from Proposition~\ref{idealcancel}, since wherever it appears, $M = \ann_R(I)$.
\end{proof}

Corollary~\ref{succinctcentral} provides a characterization of 
central minimal extensions of an arbitrary prime ring which is almost identical to the characterization 
of commutative minimal extensions of commutative domains.  
Commutative domains can, however, have non-commutative minimal
extensions which are not of the flavors appearing in Corollary~\ref{succinctcentral}, as we will see in Example~\ref{bergmanfield}.    

Our next goal is to produce examples of the extensions described in Theorem~\ref{primeext} that cannot 
be central extensions, namely those of types \ttb and \tc, and also non-central extensions of type \td.  We begin with 
an example of such a type-\td \, extension, which was brought to our attention by Alex Diesl.  By Proposition~\ref{centralprime}, we 
seek a ring with a subdirectly irreducible prime ideal which is not maximal.  

\begin{example}\cite[Exercises 3.14-3.16]{fc}   
Let $k$ be a field, and let $V$ be a countably infinite-dimensional $k$-vector space.  
The ring $R = \End_k(V)$ is a prime ring, with exactly three ideals, $0$, $R$, and the ideal $I$ consisting of all 
endomorphisms of finite rank.  Note that both $0$ and $I$ are
subdirectly irreducible prime ideals, whereas $\Max(R) = \{I\}$.  By Proposition~\ref{centralprime}, 
$E(R,I)$ is a type-\ttd minimal extension which is not central, since $\ann_R(I) = 0$.   
\end{example}

Now, note that if $S$ is a non-central type-\ttd minimal extension of a ring $R$, then $S$ 
cannot be $R$-isomorphic to the $R$-ring $R \times R/M$ for any maximal ideal $M$ (with $R$ embedded via the diagonal embedding), 
since such an extension must be central.  More generally, every central idempotent of $S$ must lie in $R$, so $S$ cannot be 
a nontrivial direct product of any two rings over $R$.  

Next, we examine extensions of type \ttb and \tc.  
Type-\ttc extensions can be used to produce type-\ttb extensions, and vice versa.  
First, let us show that any type-\ttc extension can be used to produce a type-\ttb extension.  
Indeed, suppose that $E(R,I)$ is a minimal extension with $I^2 \ne 0$,
$\ann_R(I) \ne 0$, and $\Hom_R(I,R/\ann_R(I)) = 0$.  
Then $E(R/\ann_R(I),I)$ is a minimal extension of the prime ring $S=R/\ann_R(I)$, for which $I^2 \ne 0$, $\ann_S(I) = 0$, and $\Hom_S(I,S) = 0$.  That is, the quotient $E(R/\ann_R(I),I)$ is a type-\ttb extension of the prime ring $R/\ann_R(I)$.  

Conversely, given a type-\ttb extension, one can use it to produce a type-\ttc extension as follows.  
Indeed, suppose that $E(R,I)$ is a minimal extension of $R$ with $I^2 \ne 0$, $\ann_R(I) = 0$, and $\Hom_R(I,R) = 0$.  
Find a prime ring $S$ with a noninjective surjective ring homomorphism $\varphi: S
\rightarrow R$ (e.g., take $S$ to be a suitable free ring).  Using the
$R$-rng structure of $I$, we make $I$ into an $S$-rng (which
is minimal) by defining $s \cdot i = \varphi(s) i$ and $i \cdot s = i \varphi(s)$ for $i \in I$, $s \in S$.  Then, $\ann_S(I) = \ker(\varphi)$, clearly $I^2 \ne 0$, and $\Hom_S(I,S/\ann_S(I)) \cong \Hom_S(I,R) = 0$.  That is, $E(S,I)$ is a type-\ttc minimal extension of $S$.     

We postpone producing examples of extensions of type \ttb (as well as
examples of extensions of type \tc) until
Examples~\ref{bergmanmatrix} and \ref{bergmanfield} below. 

\section{Minimal extensions of simple rings} \label{simplesection} 

In this section, we will examine minimal extensions of simple rings in more depth, producing examples of 
such extensions.  

\begin{thm} \label{simplechar}  
Let $R$ be a simple ring.  Then, up to $R$-isomorphism, every minimal
extension of $R$ must be of exactly one of the following four forms.  
\begin{enumerate}
\item[\ta] A simple ring that is a minimal extension of $R$. 
\item[\tb] $E(R,I)$, for some minimal $R$-rng $I$ such that $I^2 \ne 0$ and $I$ is not $R$-isomorphic to $R$.  
\item[\td] $R \times R$.
\item[\te] $R \propto M$, for some simple $(R,R)$-bimodule $M$.
\end{enumerate}
In addition,
\begin{enumerate}
	\item As in Theorem~\ref{primeext}, type-\tta and type-\ttb
	extensions are prime; type-\ttd extensions are semiprime but
	not prime, and type-\tte extensions are not semiprime.  
\item  Every type-\ttd extension is a central extension, and no type-\ttb
  extension is a central extension.
For type-\tte extensions, $M$ is isomorphic to $R$, as an
  $(R,R)$-bimodule, if and only if $S$ is a central extension. 
\item  If $S$ is a type-\ttb extension, $S$ has a unique
  proper nonzero ideal.
\end{enumerate}
\end{thm}
\begin{proof}
The characterization follows quickly from Theorem~\ref{primeext}.  
Let $S$ be a type-\tta extension of $R$ in the sense of
Theorem~\ref{primeext}, so any nonzero ideal of $S$ intersects $R$ nontrivially.  
If $I$ is a proper nonzero ideal of $S$, then $R \cap I$ is a 
proper nonzero ideal of $R$, which must be zero since $R$ is simple.  
We conclude that $S$ has no proper nonzero ideals and hence $S$ must be simple.  

All minimal extensions of $R$ which are not of type \tta are ideal
extensions $E(R,I)$.  
Note that $\ann_R(I)$ is a prime ideal of $R$, and hence $\ann_R(I) =
0$, since $R$ is simple.  In
particular, $R$ cannot have type-\ttc extensions.  
Type-\ttb and
type-\tte extensions of $R$, in the sense of Theorem~\ref{primeext},
reduce to those appearing in \ttb and \te, respectively.  For $R$
simple, it is clear that $0$ is the only subdirectly irreducible prime
ideal.  Thus, the only type-\ttd extension is $E(R,R)$, which is $R$-isomorphic to
$R \times R$ (the argument is similar to that in the proof of Proposition~\ref{centralprime}).  

The statements about (semi)primeness follow from those found in Theorem~\ref{primeext}, and those about centrality follow 
immediately from Proposition~\ref{centralprime} and Proposition~\ref{centralchar}.   
Finally, the statement about the numbers of ideals of $S$ follows
immediately from Proposition~\ref{idealdescription}.
\end{proof} 

\begin{remark}  A minimal extension of a simple ring can have at most
  two proper nontrivial ideals (for instance, by
  Proposition~\ref{idealdescription}, type-\ttd extensions have exactly two such
  ideals, and type-\te \, extensions have one such ideal).  
\end{remark} 

\begin{remark}  \label{primeartinian} Note that any (left or right) artinian prime ring is simple (e.g., this follows easily from \cite[Theorem 10.24]{fc}), so each type-\ttb extension of a simple ring must be non-artinian.  
\end{remark}  

Let us now construct an example, due to George Bergman, 
of a non-simple prime minimal extension of a simple ring.  
 
\begin{example} \label{bergmanmatrix}
Consider the ring homomorphism $f_n: \M_{2^n}(k) \rightarrow
\M_{2^{n+1}}(k)$ which inflates the entry $c \in k$ to 
the corresponding $2 \times 2$ scalar matrix $\begin{pmatrix} c & 0 \\ 0 &
  c \end{pmatrix}$.  Let $R$ be the direct limit of the directed system 
$$\M_2(k) \rightarrow \M_{2^2}(k) \rightarrow \M_{2^3}(k) \rightarrow
\cdots,$$ with transition maps $f_n$, as defined above.
As a direct limit of simple rings, $R$ is simple.  

Now, we will consider a slightly different direct limit; the rings
will be the same, however, the maps will no longer be ring
homomorphisms.  
Consider the idempotent $E_n = \sum_{j \text{ odd}} e_{jj} \in
\M_{2^n}(k)$, where $e_{ij}$ denotes the matrix unit with a
$1$ in position $(i,j)$, and $0$ everywhere else.  
Define the map $g_n: \M_{2^n}(k) \rightarrow \M_{2^{n+1}}(k)$ by 
$g_n(A) = E_n f_n(A) E_n$ ($g_n$ inflates each scalar
entry $c \in k$ to the $2 \times 2$ matrix $\begin{pmatrix} c & 0 \\ 0
  & 0 \end{pmatrix}$).  It is easy to verify that $f_n(A) E_n = E_n
f_n(A) = E_n f_n(A) E_n$ (check for $A = e_{rs}$ and extend linearly).  
It follows easily that $g_n(AB) = f_n(A) g_n(B) = g_n(A) g_n(B) =
g_n(A) f_n(B)$ for any $A,B \in \M_{2^n}(k)$ (with the multiplication
taking place in $\M_{2^{n+1}}(k)$).  Let $I$ be the direct limit of the
sequence of rngs with transition maps $g_n$ (which are rng
homomorphisms).  The rng $I$ can be endowed with the structure of an
$(R,R)$-bimodule as follows.    
Given $A \in R$ and $B \in I$, 
$A \in \M_{2^n}(k)$ for some $n$ and $B \in \M_{2^m}(k)$ for
some $m$.  Applying transition maps to $A$ and $B$, we may assume that
$n=m$, and then we set $A \cdot B = AB$ and $B \cdot A = BA$, where
the multiplication takes place in
$\M_{2^n}(k)$.  Note that this is well defined, since 
$f_n(A)g_n(B) = g_n(AB)$ and $g_n(B)f_n(A) = g_n(BA)$.  Using this
$(R,R)$-bimodule structure, we make $I$ into an $R$-rng.  Moreover,
$I$ is clearly simple (essentially, it is the direct limit of simple
bimodules) as an $(R,R)$-bimodule, and is hence a minimal $R$-rng.  
Finally, $I^2 \ne 0$ (nonzero idempotents abound) and it
is easy to check that $I$ is not $R$-isomorphic to $R$, as $I$ can be
shown to lack an identity element.  Indeed, if $z \in I$ were an
identity element for $I$, then $z \in \M_{2^n}(k)$ for some $n$, but
then $z$ cannot act as the identity element when viewed as $g_n(z) \in
\M_{2^{n+1}}(k)$, since $g_n(z)$ annihilates $1-E_{n+1}$.    
It follows that $E(R,I)$ is a minimal
extension of type \tb.  
\end{example}

Now let us specialize our discussion even further, to simple commutative rings, which are, of course, just 
fields.  By \cite[Lemme 1.2]{F&O}, for a field $k$, up to
$k$-isomorphism, the 
{\emph{commutative}} minimal extensions of $k$ are:  minimal field extensions (which is the commutative situation of Theorem~\ref{simplechar}), $k \times k$, and $k[x]/(x^2)$.
In light of Theorem~\ref{simplechar}, removing the commutativity
hypothesis does not change the types of rings appearing in type-\ttd or type-\tte 
dramatically, but the behavior of prime minimal extensions can be
significantly different in the noncommutative case, and we will study
this briefly.  

To start, observe that certain fields have no noncommutative minimal extensions.  
For instance, any minimal extension of a prime field 
($\Q$, or $\F_p$ for a prime $p$) is a central extension (since the subring generated by $1$ is contained in the center of the extension), and is hence commutative.

\begin{q}  If $k$ is a field with no noncommutative (or, no noncommutative prime) minimal extensions, must $k$ be a prime field?   
\end{q}

Generalizing the above, any field which is finite-dimensional over its prime field 
has the property that each of its minimal extensions is finite-dimensional over its prime field as well.  
This follows from results in \cite{Klein} and \cite{Laffey} when the
field is finite, and a result in \cite{MaxAlg}, in the case the prime
field is $\QQ$.  

Other fields can have division rings as minimal extensions.  
For instance, the division ring $\HH$ of real quaternions 
is a minimal extension of $\CC$.  
More generally, cyclic algebras (e.g., see \cite{fc}) can be used to
construct minimal extensions which are division rings.      

While on the topic of $\CC$, it is worth noting that a field is
algebraically closed if and only if each of its prime minimal
extensions if noncommutative.  Indeed, by \cite[Lemme 1.2]{F&O}, any
commutative prime minimal extension of a field $k$ is a field
extension.  But, a field extension $F$ cannot be
minimal over $k$ if $[F:k] = \infty$.  For instance, if $F$ is
algebraic over $k$, there will be a
finite-degree extension of $k$ inside $F$; if $F$ is transcendental
over $k$, and $x \in F \setminus k$, the subfield $k(x^2)$ is properly
contained in $F$.  The observation is now immediate.  

The following lemma characterizes those fields which have a simple artinian minimal extension
whose corresponding division ring is centrally finite.  

\begin{lemma} \label{finiteindex} A field $k$ has a simple artinian minimal extension of the form
  $\M_n(D)$, with $D$ a centrally finite division ring and $n>1$, if
  and only if
  $k$ has a proper subfield of finite index.  
\end{lemma}
\begin{proof}
For the forward implication, note that $k$ is clearly centralized by
$Z(D)$, and hence the subring that $k$ and $Z(D)$ generate is a
commutative subring of $\M_n(D)$.  Since $\M_n(D)$ is a minimal
extension of $k$ and $n>1$, $k$ must contain $Z(D)$.  
Note that $k \ne Z(D)$, since $Z(D) \subsetneq \T_n(Z(D))
\subsetneq \M_n(D)$.
Finally, since $\M_n(D)$, and hence $k$, is finite-dimensional over
$Z(D)$, it follows that $Z(D)$ is a subfield of $k$ having
finite index. 

For the reverse implication, suppose that $k$ has a proper subfield
$F$ of finite index, say $[k:F] = n$.  View $k$ as an $n$-dimensional
$F$-vector space, and embed $k$ in $R = \End_F(k) \cong \M_n(F)$ via
the left regular action of $k$ on itself (where $\End_F(k)$ denotes
the ring of the $F$-vector space $k$).  Now, let $S$ be a
subring of $R$ minimal among subrings of $R$ properly containing $k$.      
First, we claim that $S$ is not commutative.  To see this, note that
$[R:F] = [k:F] \cdot [\C_R(k):F]$ (see \cite[p. 105]{herstein}), 
but we already know that $[k:F] = n$ and $[R:F] = n^2$, 
so $[\C_R(k):F] = n$.  Since $k$ is
commutative, we know that $k \subseteq \C_R(k)$, and so by
comparing dimensions, we conclude that $\C_R(k) = k$.  It follows
that $S$ does not centralize $k$, and hence $S$ is not commutative.  

We claim that $S$ is prime.  By Theorem~\ref{simplechar}, the only 
semiprime non-prime minimal extensions of $k$ are $k$-isomorphic to $k
\times k$, 
which is commutative, so we need only show that $S$ is semiprime.  
To see this, let $M$ be any $(k,k)$-subbimodule of $R = \End_F(k)$,
and let $0 \ne \varphi \in M$.
Viewing $\varphi$ as an element of $R$,   
we may thus find $x \in k$ such that $\varphi(x) = y \ne 0$.  Let $\alpha$ denote left multiplication 
by $xy^{-1}$, as an element of the copy of $k$ inside $R$.  
Then, $\varphi \alpha \varphi \ne 0$, since $\varphi(\alpha(\varphi(x))) = y \ne 0$.  
Since $\varphi$ and  $\alpha \varphi$ are both elements of $M$, we conclude that $M^2 = 0$ only if $M = 0$.  
It follows that $S$ must be semiprime, and hence prime.  

Next, we claim that $S$ is not a division ring.  Indeed, note that the
$F$-subspace $W$ of $R$ consisting of matrices with no nonzero entries in the first row 
has $F$-dimension $n(n-1)$ and consists entirely of zero divisors.  
%
On
the other hand, $S$ is an $F$-subspace with $\dim_F(S) \ge 2n$ (since
$\dim_k(S) > 1$ and $\dim_F(k) = n$).  Since $\dim_F(R) = n^2 <
\dim_F(S) + \dim_F(W)$, we conclude that $S$ intersects $W$
nontrivially, and hence $S$ is not a division ring.     

Finally, note that $S$ is finite-dimensional over $F$, and hence $S$ is clearly artinian.  
Since $S$ is prime, we conclude from Remark~\ref{primeartinian} that
$S$ is simple artinian.  
Thus,  
$S \cong \M_m(D)$ for some division ring $D$ containing $F$.
Moreover, $D$ must be centrally finite, since $R$ is 
finite-dimensional over $F$.   
Also note that by the Double Centralizer
Theorem (e.g., see \cite[p.105]{herstein}), the centralizer of $S$ in $R$
is a maximal subfield of $k$.    
\end{proof}
 
\begin{remark}  
While the proof of Lemma~\ref{finiteindex} shows that every simple
artinian minimal extension of a field $k$ is attached to a maximal
finite index subfield of $k$, this correspondence is not unique.  
The proof above shows, in fact, that if $F \subseteq k$ is a maximal subfield of
finite index $n$, then $\End_F(k) \cong \M_n(F)$ is a minimal
extension of $k$.  Were this correspondence unique, we would have been able to sharpen the statement of Lemma~\ref{finiteindex} by replacing  ``centrally finite division
ring'' with ``field''.  This is not the case, however.  For instance, let $\HH$ denote
the division ring of ordinary quaternions with rational coefficients.
There exist subfields of $\M_2(\HH)$ which are degree-four extensions
of $\QQ$, and in which $\QQ$ is a maximal subfield; for instance, the
subring generated by $\begin{pmatrix} -3 -i + 2j + 3k & 3+3i+3j-k \\ 3
  - 2i + 3j - k & -2i -3j + k \end{pmatrix}$ is such a field (the
associated Galois group is $S_4$, and all of its subgroups of index
four are maximal).         
\end{remark}

\begin{q}   Does there exist a field with a minimal ring extension which is 
a centrally infinite division ring?  
\end{q}  

\begin{q}Does there exist a field with a minimal ring extension of the form $\M_n(D)$, where $n>1$ and $D$ is a centrally infinite division ring? 
\end{q}   

\begin{q}Does there exist a field with a minimal ring extension which is a non-artinian simple ring?
\end{q}  

The following example due to George Bergman shows, 
in particular, that any field which is purely transcendental over a 
subfield has a non-simple prime minimal extension.  

\begin{example}  \label{bergmanfield}
Let $F$ be any field and let $F((t)) = F[[t]][t^{-1}]$ denote the field of formal Laurent
series in one variable $t$ over $F$.  Let $\tr: F((t)) \rightarrow F$
be the $F$-linear map sending a formal Laurent series to its constant
coefficient.  Also, the degree of a Laurent series $\sum_{i \in \Z}
a_i t^i$ is the least $i$ for which $a_i \ne 0$, and is $-\infty$
if no such integer exists (which only happens for the series $0$).     

Let $k$ be any subfield of
$F((t))$ which contains $F(t)$.  
Consider the rng $I$ whose additive group is that of $k \otimes_F k$,
but with multiplication defined by $(a \otimes b)(c \otimes d) =
\tr(bc) (a \otimes d)$.  
It is straightforward to check that this, together with the left and
right $k$-actions $a(b \otimes c) = (ab) \otimes c$ and $(a \otimes
b)c = a \otimes (bc)$ endows $I$ with the structure of an $k$-rng.  
Note that $1 \otimes 1$ is a nonzero idempotent of $I$, so $I^2 \ne
0$.  Also note that $t$ acts noncentrally, since $t(1 \otimes 1) = t
\otimes 1 \ne 1 \otimes t = (1 \otimes 1)t$.  

We claim that $I$ has no nonzero $k$-subrngs.  Indeed, we will show
that for any pair of nonzero $x,y \in I$, there exist $a,b,c \in k$
for which $axbyc = 1 \otimes 1$.  From this, it follows that the
$k$-subrng of $I$ generated by any nonzero $x \in I$ contains $1 \otimes
1$, hence contains all tensors, and hence contains $I$.  
To prove the claim, let $x,y \in I$ be nonzero.  We may write 
$x = \sum_{i=1}^n f_i \otimes g_i$ and $y = \sum_{j=1}^m r_j \otimes
s_j$ where each $f_i, g_i, r_j, s_j \in k$ for each $i,j$.  
Any finite-dimensional $F$-subspace of $k$ has a basis consisting of
series which have distinct degrees.  Using this and $F$-linearity, we may
assume that $\deg(g_i) \ne \deg(g_{i'})$ if $i \ne i'$, and $\deg(r_j)
\ne \deg(r_{j'})$ if $j \ne j'$.  Moreover, discarding any nonzero
terms, we may assume that each $f_i,g_i,r_j$ and $s_j$ is nonzero.  We
may further assume that $\deg(g_1) < \deg(g_i)$ for each $i>1$ and
$\deg(r_1) < \deg(r_j)$ for each $j>1$.  Set $a = f_1^{-1}$, $b =
g_1^{-1} r_1^{-1}$, and $c = s_1^{-1}$.  
Now, $xby = \sum_{i,j} \tr(g_i b r_j) f_i \otimes s_j$.  Note that
$g_i b r_j$ has degree 
$$\deg(g_i) + \deg(b) + \deg(r_j) = \deg(g_i) - \deg(g_1) - \deg(r_1) +
\deg(r_j),$$ which is strictly positive unless $i=j=1$.  Thus, $xby =
\tr(g_1 b r_1) f_1 \otimes s_1 = f_1 \otimes s_1$.  It follows that
$axbyc = 1 \otimes 1$.  We conclude that $I$ has no proper nonzero $k$-subrngs.
Finally, since $t$ acts noncentrally on $I$, $I$ cannot be $k$-isomorphic to
$k$.  Thus, $E(k,I)$ is a minimal extension of $k$ of type \tb.          
\end{example}

\begin{example}  Given a field $k$, the commutative domain $k(t)[x]$ has a
  type-\ttc minimal extension, using Example~\ref{bergmanfield}
  together with the argument at the end of
  Section~\ref{primesection}.  
\end{example}

We close with some questions regarding non-simple prime minimal
extensions of fields.  
\begin{q}  \label{whichnonsimple} Which fields possess a non-simple prime minimal extension?  Are
  these precisely the fields of transcendence degree at least 1?  
\end{q}

\begin{q}  \label{passtoalg}  If $k$ is a field which has a non-simple prime minimal extension, and $F$ is an 
algebraic extension of $k$, must $F$ possess a non-simple prime minimal extension?  
\end{q}  
Since Example~\ref{bergmanfield} shows that any field $k$ which is purely transcendental over a subfield has a non-simple 
prime minimal extension, an affirmative answer to Question~\ref{passtoalg} would show that any field of 
positive transcendence degree has a non-simple prime minimal
extension.  To answer Question~\ref{passtoalg} affirmatively, it would
suffice to show that if $k$ is a field with a non-simple minimal
extension $I$, and $F = k(a)$ is a minimal finite field extension of
$k$ with a primitive element $a \in F$, then $F$ has a non-simple
minimal extension $I'$, together with a $k$-rng injection $I
\rightarrow I'$ (one can then use transfinite induction to complete
the proof).

\bibliography{minimal}

\noindent
Center for Communications Research\\
4320 Westerra Court\\
San Diego, CA 92126-1967\\
USA\\
Email: {\tt dorsey@ccrwest.org} \\

\noindent
Department of Mathematics \\
University of Southern California \\
Los Angeles, CA 90089 \\
USA \\
Email: {\tt mesyan@usc.edu}

\end{document}